\documentclass[11]{article}

\usepackage{amsmath}
\usepackage{amsfonts}
\usepackage{amsthm}

\usepackage{fullpage}

\usepackage{hyperref}

\newcommand{\Tr}{\ensuremath{\mathrm{Tr}}}
\newcommand{\Airy}{\ensuremath{\mathrm{Ai}}}
\newcommand{\diag}{\ensuremath{\textrm{diag}}}
\newcommand{\PP}{\mathbb{P}}
\newcommand{\PPn}{\mathbb{P}_n}

\newcommand{\PPDyson}{\mathbb{P}^{\mathrm{Dyson}}_n}
\newcommand{\PPAiry}{\mathbb{P}^{\mathrm{Airy}}}

\newtheorem{thm}{Theorem}
\newtheorem{cor}{Corollary}
\newtheorem{lemma}{Lemma}

\theoremstyle{remark}
\newtheorem{remark}{Remark}

\title{A PDE for the Multi-Time Joint Probability of the Airy Process}
\author{Dong Wang \footnote{Department of Mathematics, Brandeis University, Waltham, MA 02454, USA. \href{mailto:wangdong@brandeis.edu}{wangdong@brandeis.edu}}}
\begin{document}

\maketitle

\begin{abstract}
This paper gives a PDE for multi-time joint probability of the Airy
process, which generalizes Adler and van Moerbeke's result on the
$2$-time case. As an intermediate step, the PDE for the multi-time
joint probability of the Dyson Brownian motion is also given.
\end{abstract}

\section{Introduction}

The Airy process can be defined purely stochastically as the limit
of the Dyson Brownian motion, as we are going to do later. However,
it also appears in various statistical physical models, such as the
polynuclear growth process \cite{Prahofer-Spohn02},
\cite{Johansson03} and the Domino tiling model \cite{Johansson05}.
Since the Airy process is stationary with continuous sample path
\cite{Prahofer-Spohn02}, we can pick any time $t$ and consider the
probability of all particles being in $(-\infty, u)$, denoted by
$\PP(u)$, and find that the probability is given by the GUE
Tracy-Widom distribution \cite{Tracy-Widom94}
\begin{equation}
\PP(u) = e^{-\int^{\infty}_u (s-u)q^2(s)ds},
\end{equation}
with $q(s)$ the solution of the Painlev\'{e} II equation
\begin{equation}
q''(s) = sq(s)+2q^2(s), \qquad q(s) \simeq
\begin{cases}
\frac{e^{-(2/3)s^{3/2}}} {2\sqrt{\pi}s^{1/4}} & \textrm{for $s
\rightarrow \infty$,} \\
\sqrt{-s/2} & \textrm{for $s \rightarrow -\infty$.}
\end{cases}
\end{equation}

In their study of the joint probability for several times of the
Airy process, Pr\"{a}hofer and Spohn \cite{Prahofer-Spohn02} posed
the problem to find a PDE for the joint probability. Adler and van
Moerbeke \cite{Adler-van_Moerbeke05} solved the problem for the
$2$-time case, and assuming a plausible conjecture of the boundary
condition, got the asymptotic expansion of the probability function
$\PP(t,u, v)$, which is the probability that all particles are in
$(-\infty, u)$ initially and in $(-\infty, v)$ after a large time
$t$. Their solution was obtained by a previous result of them on the
spectrum of coupled random matrices \cite{Adler-van_Moerbeke99}.
They regarded the joint distribution for the Dyson Brownian motion
of $2$-time as a $\tau$ function of the two-Toda lattice, and
construct a PDE with variables in times and boundary points of the
Dyson Brownian motion as a consequence of identities for $\tau$
functions and Virasoro identities specific to the situation. Then
they got the PDE for the Airy process by taking the limit.

This paper generalizes their result to the multi-time case, and the
technical heart is the same identity for $\tau$ functions, although
in the generalized case we need more elaborate work to fit
differential operators in times and boundary points of the Dyson
Brownian motion into the structure of two-Toda $\tau$ functions.

After the description of the problem, we state the PDEs for both the
Dyson Brownian motion with finite particles and its limit, the Airy
process with infinitely many particles, and an example for the
$3$-time ($m=2$) case for the Airy process. Section
\ref{PDE_for_the_finite_Dyson_process} derives the result for the
Dyson process and section \ref{PDE_for_the_Airy_process} derives the
result for the Airy process by taking limit.

\subsection{Description of the model}
\label{Description_of_the_model}

The free Brownian motion process is determined by the transition
probability distribution
\begin{equation}
P(t,\bar{X},X) = \frac{1}{\sqrt{(2\pi
t)/\beta}}e^{-\frac{(X-\bar{X})^{2}}{2t/\beta}},
\end{equation}
where $\bar{X}$ and $X$ are initial and terminal coordinates of the
particle, and $\beta$ is the diffusion constant. The probability
distribution $P(t,\bar{X},X)$, as a function of $t$ and $X$,
satisfied the diffusion equation
\begin{equation}
\frac{\partial P}{\partial t} = \frac{1}{2\beta}
\frac{\partial^{2}}{\partial X^{2}}P.
\end{equation}
If we add a harmonic potential $\rho X^2/2$ to the process, then the
probability distribution $P(t,\bar{X},X)$ satisfies (see e.g.
\cite{Feller68})
\begin{equation}
\frac{\partial P}{\partial t} = \left( \frac{1}{2\beta}
\frac{\partial^{2}}{\partial X^{2}} - \frac{\partial}{\partial X} (-
\rho X)\right)P,
\end{equation}
and the process is determined by ($c = e^{-\rho t}$)
\begin{equation}
P(t,\bar{X},X) = \frac{1}{\sqrt{\frac{\pi (1-c^{2})}{\rho \beta}}}
e^{-\frac{(X-c\bar{X})^{2}}{(1-c^{2})/\rho \beta}}.
\end{equation}
While the free Brownian motion process is dispersive, the Brownian
motion process in the harmonic potential well has a stationary
distribution
\begin{equation}
P(X) = \frac{e^{-\rho \beta X^{2}}}{\sqrt{\pi/\rho \beta}}.
\end{equation}

Now we can define the Ornstein-Uhlenbeck process
\cite{Ornstein-Uhlenbeck30} of an $n \times n$ Hermitian matrix $B$,
in which all the $n^2$ real variables---$n$ for real diagonal
entries, $n(n-1)/2$ for the real parts of off diagonal entries, and
the other $n(n-1)/2$ for the imaginary parts of them---are in
independent Brownian motion in harmonic potential wells. The $\rho$
for them are uniformly $1$, and $\beta$ are $1$ for the $n$ diagonal
variables and $2$ for the $n(n-1)$ off diagonal variables. Therefore
for $i$, $j$ in $\{ 1, \dots, n \}$, ($c = e^{-t}$)
\begin{equation}
\begin{cases}
P_{ii}(t,\bar{B}_{ii},B_{ii}) = \frac{1}{\sqrt{\pi (1-c^{2})}}
e^{-\frac{(B_{ii}-c\bar{B}_{ii})^{2}}{1-c^{2}}}, \\
P_{ij\Re}(t,\Re \bar{B}_{ij},\Re B_{ij}) = \frac{1}{\sqrt{\pi
(1-c^{2})/2}}
e^{-\frac{(\Re B_{ii}-c\Re \bar{B}_{ii})^{2}}{(1-c^{2})/2}}, \\
P_{ij\Im}(t,\Im\bar{B}_{ij},\Im B_{ij}) = \frac{1}{\sqrt{\pi
(1-c^{2})/2}} e^{-\frac{(\Im B_{ii}-c\Im
\bar{B}_{ii})^{2}}{(1-c^{2})/2}},
\end{cases}
\end{equation}
and we can write the joint transition probability distribution as
\begin{equation}
P(t,\bar{B},B) = \prod^n_{i=1}P_{ii} \prod_{1 \leq i < j \leq n}
\left( P_{ij\Re}P_{ij\Im} \right) =
\frac{C^{-1}}{(1-c^{2})^{n^{2}/2}} e^{-\frac{\Tr
(B-c\bar{B})^{2}}{1-c^{2}}}.
\end{equation}

We consider the multi-time transition function with the initial
state $B_0$ at $t_0 = 0$, the terminal state $B_m$ and a series of
intermediate states $B_1, \dots, B_{m-1}$, and the time between
state $B_0$ and $B_i$ being $t_i$, if we denote
\begin{equation}
s_i =
\begin{cases}
0 & i=0, \\
t_1 & i=1, \\
t_i - t_{i-1} & i=2, \dots, m,
\end{cases}
\end{equation}
and
\begin{equation}
c_i = e^{-s_i},
\end{equation}
then
\begin{equation}
P(t_{1},\dots,t_{m},B_{0},\dots,B_{m}) = C^{-1} \prod^{m}_{i=1}
e^{-\Tr\frac{(B_{i}-c_{i}B_{i-1})^{2}}{1-c^{2}_{i}}}.
\end{equation}
The Ornstein-Uhlenbeck process has a stationary distribution
\begin{equation}
\label{eq:stationary_Ornstein-Uhlenbeck_distribution}
P(B) = C^{-1}e^{-\Tr B^2}.
\end{equation}

Since the Ornstein-Uhlenbeck process is invariant under the unitary
transiformation, we define the process of the eigenvalues as the
Dyson Brownian motion process \cite{Dyson62}, whose multi-time
transition probability distribution is ($0 = t_0 < t_1 < \dots, <
t_m$)
\begin{equation}
P(t_{1},\dots,t_{m},\lambda^{(0)},\dots,\lambda^{(m)}) =
\textrm{\begin{minipage}[c]{0.5\textwidth}{The transition
probability of the $n \times n$ Hermitian matrix with eigenvalues
initially
$\lambda^{(0)}=(\lambda^{(0)}_{1},\dots,\lambda^{(0)}_{n})$ and
$\lambda^{(1)}$ after time $t_{1}$, $\lambda^{(2)}$ after time
$t_{2}$ \dots and $\lambda^{(m)}$ after the total time $t_{m}$.}
\end{minipage}}
\end{equation}
If we change the coordinates of the $\mathbb{R}^{n^2}$ space of $n
\times n$ Hermitian matrices in to the eigenvalue-angle coordinates
$\lambda_1, \dots, \lambda_n$, $\theta_1, \dots, \theta_{n(n-1)}$,
with the Jacobian identity (see e.g. \cite{Mehta04})
\begin{equation}
\prod^{n}_{i=1}dx_{ii} \prod_{1 \leq i < j \leq n} \left(
d\Re(x_{ij}) d\Im(x_{ij}) \right) = {V(\lambda)}^{2}
\prod^{n}_{i=1}d\lambda_{i}\prod^{n(n-1)}_{i=1}d\theta_{i},
\end{equation}
where $V(\lambda) = \prod_{1 \leq i < j \leq n} (\lambda_i -
\lambda_j)$ is the Vandermonde, we find the explicit formula for
$P(t_{1},\dots,t_{m},\lambda^{(0)},\dots,\lambda^{(m)})$:
\begin{equation}
P(t_{1},\dots,t_{m},\lambda^{(0)},\dots,\lambda^{(m)}) = \frac{1}{C}
\int\dots\int \prod^{m}_{i=1} e^{-\Tr \frac{\left(
B(\lambda^{(i)},\theta^{(i)}) -
c_{i}B(\lambda^{(i-1)},\theta^{(i-1)}) \right)^{2}}{1-c^{2}_{i}}}
\prod^{m}_{i=1}V(\lambda^{(i)})^{2} \prod^{m}_{i=1}d\theta^{(i)},
\end{equation}
where $\theta^{(0)}$ does not appear in the integral since the
transition probability is independent of $\theta^{(0)}$ for the
unitary invariant property.

By the Harish-Chandra-Itzykson-Zuber (HCIZ) formula
\cite{Bessis-Itzykson-Zuber80}
\begin{equation}
\int_{U(n)}e^{\Tr (XUYU^{-1})} dU =
C\frac{\det(e^{x_{i}y_{j}})}{V(x)V(y)},
\end{equation}
where $X=\diag(x_{1},\dots,x_{n})$ and $Y=\diag(y_{1},\dots,y_{n})$
are diagonal matrices, we can evaluate the multi-time transition
probability density as
\begin{multline}
P(t_{1},\dots,t_{m},\lambda^{(0)},\dots,\lambda^{(m)}) = \frac{1}{C}
V(\lambda^{(0)})^{-1} V(\lambda^{(m)}) \prod^{m}_{l=1} \det \left(
e^{\frac{2c_{l}}{1-c^{2}_{l}} \lambda^{(l-1)}_i \lambda^{(l)}_j} \right)\\
e^{-\frac{c^{2}_{1}}{1-c^{2}_{1}}
\sum^{n}_{i=1}{\lambda^{(0)}_{i}}^{2}} \prod^{m-1}_{l=1} \left(
e^{-\left( \frac{1}{1-c^{2}_{l}} + \frac{c^{2}_{l+1}}{1-c^{2}_{l+1}}
\right) \sum^{n}_{i=1} {\lambda^{(l)}_{i}}^{2}} \right)
e^{-\frac{1}{1-c^{2}_{m}} \sum^{n}_{i=1} {\lambda^{(m)}_{i}}^{2}}.
\end{multline}

If we take the initial state with eigenvalue $\lambda^{(0)}$ from
the stationary distribution
(\ref{eq:stationary_Ornstein-Uhlenbeck_distribution}), which is
\begin{equation}
\tilde{P}(\lambda^{(0)}) =\frac{1}{C} V(\lambda^{(0)})^{2}
e^{-\sum^{n}_{i=1} {\lambda^{(0)}_{i}}^{2}},
\end{equation}
We get the multi-time correlation function in the stationary Dyson
process
\begin{equation}
\begin{split}
\tilde{P}(t_{1},\dots,t_{m},\lambda^{(0)},\dots,\lambda^{(m)}) = &
\tilde{P}(\lambda^{(0)})
P(t_{1},\dots,t_{m},\lambda^{(0)},\dots,\lambda^{(m)}) \\
= & \frac{1}{C} V(\lambda^{(0)}) V(\lambda^{(m)}) \prod^{m}_{l=1}
\det \left( e^{\frac{2c_{l}}{1-c^{2}_{l}} \lambda^{(l-1)}_i
\lambda^{(l)}_j} \right) \\
& e^{-\frac{1}{1-c^{2}_{1}} \sum^{n}_{i=1} {\lambda^{(0)}_{i}}^{2}}
\prod^{m-1}_{l=1} \left( e^{-\left( \frac{1}{1-c^{2}_{l}} +
\frac{c^{2}_{l+1}}{1-c^{2}_{l+1}} \right) \sum^{n}_{i=1}
{\lambda^{(l)}_{i}}^{2}} \right) e^{-\frac{1}{1-c^{2}_{m}}
\sum^{n}_{i=1} {\lambda^{(m)}_{i}}^{2}}.
\end{split}
\end{equation}

If we want to find the probability of all $\lambda^{(l)}_{i}$'s
being in $U^{(l)} = (a^{(l)}_{1}, a^{(l)}_{2})\cup \dots
\cup(a^{(l)}_{2r_{l}-1}, a^{(l)}_{2r_{l}})$, with $-\infty \leq
a^{(l)}_{1} < a^{(l)}_{2} < a^{(l)}_{3} < \dots < a^{(l)}_{2r_{l}}
\leq \infty$, for $l = 0,1,\dots,m $ and $i = 1,\dots,n$, which is
\begin{equation} \label{eq:definition_of_tau_n}
\PPDyson(t_1, \dots, t_m; a^{(0)}_1, \dots, a^{(m)}_{2r_m}) =
\idotsint_{{U^{(0)}}^n \times \dots \times {U^{(m)}}^n}
\tilde{P}(t_{1},\dots,t_{m},\lambda^{(0)},\dots,\lambda^{(m)})
\prod^{m}_{l=0}\prod^{n}_{k=1}d\lambda^{(l)}_{k},
\end{equation}
we can simplify it by the symmetry of $\lambda^{(l)}_{1}, \dots,
\lambda^{(l)}_{n}$, for all $l = 0, \dots, m$ and get
\begin{multline}
\PPDyson(t_i; a^{(l)}_i) = \frac{1}{C} \idotsint_{{U^{(0)}}^n \times
\dots \times {U^{(m)}}^n} V(\lambda^{(0)}) V(\lambda^{(m)}) \\
e^{-\frac{1}{1-c^{2}_{1}}\sum^{n}_{k=1}{\lambda^{(0)}_{k}}^{2}}
\prod^{m-1}_{l=1}
e^{-\frac{1-c^{2}_{l}c^{2}_{l+1}}{(1-c^{2}_{l})(1-c^{2}_{l+1})}
\sum^{n}_{k=1}{\lambda^{(l)}_{k}}^{2}}
e^{-\frac{1}{1-c^{2}_{m}}\sum^{n}_{k=1}{\lambda^{(m)}_{k}}^{2}}
\prod^{m}_{l=1}
e^{\frac{2c_{l}}{1-c^{2}_{l}}\sum^{n}_{k=1}\lambda^{(l-1)}_{k}\lambda^{(l)}_{k}}
\prod^{m}_{l=0}\prod^{n}_{k=1}d\lambda^{(l)}_{k}.
\end{multline}
We are going to give a PDE satisfied by $\log\PPDyson$ with
variables $t_i$ and $a^{(l)}_i$.

The Airy process can be defined as the limit of the Dyson process on
the edge \cite{Adler-van_Moerbeke05}. As $n \rightarrow \infty$, we
can prove that the right-most particle in the Dyson process is
almost surely around $\sqrt{2n}$ with the fluctuation scale
$n^{1/6}$ \cite{Prahofer-Spohn02}, \cite{Adler-van_Moerbeke05}. If
we take the rescaling
\begin{align}
\bar{t}_l = & n^{1/3}t_l, \label{eq:first scaled variable for Airy process}\\
\bar{\lambda}^{(l)}_k = & \sqrt{2}n^{1/6} (\lambda^{(l)}_k -
\sqrt{2n}) \\
\bar{a}^{(l)}_i = & \sqrt{2}n^{1/6} (a^{(l)}_i - \sqrt{2n}),
\label{eq:last scaled variable for Airy process}
\end{align}
then for fixed $\bar{t}_l$ and $\bar{a}^{(l)}_i$, $\PPDyson$
converges to a function defined by the Fredholm determinant of a
matrix integral operator \cite{Forrester-Nagao-Honner99},
\cite{Macedo94}, \cite{Prahofer-Spohn02}, \cite{Johansson03}
\begin{equation} \label{eq:definition of Airy process}
\lim_{n \rightarrow \infty} \left. \PPDyson \right|_{\substack{t_l =
\bar{t}_l / n^{1/3} \\ a^{(l)}_i = \sqrt{2n} + \bar{a}^{(l)}_i /
(\sqrt{2}n^{1/6})}} = \PPAiry(\bar{t}_1, \dots, \bar{t}_m,
\bar{a}^{(0)}_1, \dots, \bar{a}^{(m)}_{2r_m}) = \det \left( I -
(\chi^c_i K^A_{ij} \chi^c_j)_{1 \leq i, j \leq m} \right),
\end{equation}
where $\chi^c_l$ is an indicator function defined as
\begin{equation}
\chi^c_l (t) =
\begin{cases}
0 & t \in \bigcup^{t_l}_{i=1}
(a_{2i-1}, a_{2i}), \\
1 & \textrm{otherwise.}
\end{cases}
\end{equation}
and ($\Airy$ is the Airy function)
\begin{equation}
K^A_{ij}(x,y) =
\begin{cases}
\int^{\infty}_0 \Airy(x+z)\Airy(y+z) dz & \textrm{if $i=j$,} \\
\int^{\infty}_0 e^{-z (\bar{t}_i - \bar{t}_j)} \Airy(x+z)\Airy(y+z)
dz & \textrm{if $i>j$,} \\
-\int^0_{-\infty} e^{z (\bar{t}_j - \bar{t}_i)} \Airy(x+z)\Airy(y+z)
dz & \textrm{if $i<j$.}
\end{cases}
\end{equation}
Then we can define the Airy process, which contains infinitely many
particles by the probability function (\ref{eq:definition of Airy
process}). Furthermore, we are going to give a PDE satisfied by
$\tau$ with variables $\bar{t}_i$ and $\bar{a}^{(l)}_i$.

\begin{remark} \label{remark:on_al1}
To make the definition (\ref{eq:definition of Airy process})
meaningful, we need $\bar{a}^{(l)}_1$ to be $-\infty$ for all $l$.
Otherwise the left hand side of (\ref{eq:definition of Airy
process}) is $0$ and the right hand side is not well defined.
\end{remark}

\subsection{Statement of main results}

With notations defined in subsection
(\ref{Description_of_the_model}), we define differential operators
($l = 0, 1, \dots, m$)
\begin{equation}
D^{l,1} = \sum^{2r_{l}}_{i=1}\frac{\partial}{\partial a^{(l)}_{i}},
\qquad D^{l,2} =
\sum^{2r_{l}}_{k=1}a^{(l)}_{k}\frac{\partial}{\partial a^{(l)}_{i}},
\end{equation}
if all $a^{(l)}_i$ are finite; otherwise we drop the $a^{(l)}_1$
(resp. $a^{(l)}_{2r_l}$) part if $a^{(l)}_1 = -\infty$ (resp.
$a^{(l)}_{2r_l} = \infty$). And then denote
\begin{align}
\mathcal{A}_1 = & \sum^m_{l=0} e^{-t_l} D^{l,1}, \\
\mathcal{B}_1 = & \sum^m_{l=0} e^{t_l - t_m} D^{l,1}, \\
\mathcal{A}_2 = & \sum^m_{l=0} e^{-2t_l} D^{l,2} + \sum^m_{l=1} (1 -
e^{-2t_l}) \frac{\partial}{\partial t_l} - e^{-2t_m}, \\
\mathcal{B}_2 = & \sum^m_{l=0} e^{2(t_l - t_m)} D^{l,2} +
\sum^m_{l=1} (e^{2(t_l - t_m)} - e^{-2t_m}) \frac{\partial}{\partial
t_l} - e^{-2t_m}.
\end{align}
We now state
\begin{thm}[Dyson Brownian motion]
Given $t_1, \dots, t_m$, the logarithm of the joint distribution for
the stationary Dyson Brownian motion $\PPDyson$ defined in
(\ref{eq:definition_of_tau_n}) (abbreviated as $\log\mathbb{P}_n$)
satisfies a third order non-linear PDE in times and boundary points
of $U^{(l)}$
\begin{equation} \label{eq:compact equation for Dyson process}
\mathcal{A}_1 \frac{\mathcal{B}_2\mathcal{A}_1 \log \PPn}
{\mathcal{B}_1\mathcal{A}_1 \log \PPn +2n e^{-t_m}} = \mathcal{B}_1
\frac{\mathcal{A}_2\mathcal{B}_1 \log \PPn}
{\mathcal{A}_1\mathcal{B}_1 \log \PPn +2n e^{-t_m}}.
\end{equation}
\end{thm}

Similarly with the notations
\begin{align}
\mathcal{D} = & \sum^m_{l=0} D^{l,1}, \\
\mathcal{D}_{1L} = & \sum^m_{l=0} (t_m - t_l) D^{l,1}, \\
\mathcal{D}_{2R} = & \sum^m_{l=0} t_l D^{l,1}, \\
\mathcal{D}_1 = & = \mathcal{D}_{1L} - \mathcal{D}_{1R} = \sum^m_{l=0}
(t_m - 2t_l) D^{l,1}, \\
\mathcal{D}_2 = & \sum^m_{l=0} ((t_m - t_l)^2 + t^2_l) D^{l,1}, \\
\mathcal{D}_3 = & \sum^m_{l=0} ((t_m - t_l)^3 - t^3_l) D^{l,1}, \\
\mathcal{E} = & \sum^m_{l=0} D^{2,1}, \\
\mathcal{E}_1 = & \sum^m_{l=0} (t_m - 2t_l) D^{l,2}, \\
\mathcal{T}_1 = & 2\sum^m_{l=1} t_l \frac{\partial}{\partial t_l}, \\
\mathcal{T}_2 = & 2\sum^m_{l=1} t_l(t_m - t_l)
\frac{\partial}{\partial t_l},
\end{align}
we state the result for the Airy process
\begin{thm}[Airy process]
Given $t_1, \dots, t_m$, the logarithm of the joint distribution for
the Airy process $\PPAiry$ defined in (\ref{eq:definition of Airy
process}) (abbreviated as $\log\mathbb{P}$) satisfies a third order
non-linear PDE \footnote{in terms of the Wronskian, defined in
subsection \ref{notational_convenience}} in times and boundary
points of $U^{(l)}$
\begin{equation} \label{eq:equation_for_Airy_process}
\mathcal{D}^2 [\mathcal{E}_1 + \mathcal{D}_3 + \mathcal{T}_2]
\log\mathbb{P} - \mathcal{D}\mathcal{D}_1 [\mathcal{E} +
\mathcal{D}_2 + \mathcal{T}_1] \log\mathbb{P}
-2\mathcal{D}_{1L}\mathcal{D}_{1R}\mathcal{D}_1 \log\mathbb{P} = \{
\mathcal{D}^2 \log\mathbb{P}, \mathcal{D}\mathcal{D}_1
\log\mathbb{P} \}_\mathcal{D}.
\end{equation}
\end{thm}

In the case of $m=1$, our result agrees with that in
\cite{Adler-van_Moerbeke05}. Especially, if $U^{(0)}=(-\infty, u)$,
$U^{(1)}=(\infty, v)$ and denote $t_1=t$, then the result for
$\log\PPAiry(t,u,v)$ is
\begin{cor}[\cite{Adler-van_Moerbeke05}]
The logarithm of the $2$-time joint probability for the Airy process
$\PPAiry(t,u,v)$ (abbreviated as $\log\mathbb{P}$) satisfies a third
order non-linear PDE in variables $u$, $v$ and $t$
\begin{multline}
\Bigg[ (v-u)\left[ \frac{\partial}{\partial u} +
\frac{\partial}{\partial v} \right] \frac{\partial^2}{\partial u
\partial v} + t\left[ \frac{\partial^2}{\partial u^2} -
\frac{\partial^2}{\partial v^2} \right] \frac{\partial}{\partial t}
+ t^2\left[ \frac{\partial}{\partial u} - \frac{\partial}{\partial
v} \right] \frac{\partial^2}{\partial u \partial v} \Bigg]
\log\mathbb{P} = \\
\frac{1}{2} \left\{ \left[ \frac{\partial^2}{\partial u^2} -
\frac{\partial^2}{\partial v^2} \right]\log\mathbb{P}, \left[
\frac{\partial}{\partial u} + \frac{\partial}{\partial v} \right]^2
\log\mathbb{P} \right\}_{\frac{\partial}{\partial u} +
\frac{\partial}{\partial v}}.
\end{multline}
\end{cor}

In the $m=2$ case, if $U^{(0)}=(-\infty, u)$, $U^{(1)}=(\infty, v)$,
$U^{(2)}=(\infty, w)$, $t_1=t$ and $t_2=s$ the result for
$\log\PPAiry(t,s,u,v,w)$ is
\begin{cor}
The logarithm of the $3$-time joint probability for the Airy process
$\PPAiry(t,s,u,v.w)$ (abbreviated as $\log\mathbb{P}$) satisfies a
third order non-linear PDE in variables $u$, $v$, $w$, $t$ and $t$
($D = \frac{\partial}{\partial u} + \frac{\partial}{\partial v} +
\frac{\partial}{\partial u}$)
\begin{multline}
\Bigg[ t(u-v)\frac{\partial^2}{\partial u \partial v} + s(u-w)
\frac{\partial^2}{\partial u \partial w} +
(s-t)(v-w)\frac{\partial^2}{\partial u \partial w} + \left[
-s\frac{\partial}{\partial u} + (2t-s)\frac{\partial}{\partial v} +
s\frac{\partial}{\partial w} \right] \left[
t\frac{\partial}{\partial t} + s\frac{\partial}{\partial s} \right]
\\
+ t(s-t)D\frac{\partial}{\partial t} \Bigg]D \log\mathbb{P} + \left[
-t^3\frac{\partial^3}{\partial u^2 \partial v} -
s^3\frac{\partial^3}{\partial u^2 \partial w} +
t^3\frac{\partial^3}{\partial u \partial v^2} \right. \\
\left. + (2t-s)(2s-t)(s+t)\frac{\partial^3}{\partial u \partial v
\partial w} + s^3\frac{\partial^3}{\partial u \partial w^2} -
(s-t)^3\frac{\partial^3}{\partial v^2 \partial w} +
(s-t)^3\frac{\partial^3}{\partial v \partial w^2} \right]
\log\mathbb{P} = \\
\frac{1}{2} \left\{ \left[ -s\frac{\partial}{\partial u} +
(2t-s)\frac{\partial}{\partial v} + s\frac{\partial}{\partial w}
\right]D \log\mathbb{P}, D^2 \log\mathbb{P} \right\}_D.
\end{multline}
\end{cor}

\subsection{Notational convenience} \label{notational_convenience}
Throughout this paper, parentheses $(\dots)$ always include numbers
and functions; brackets $[\dots]$ always include operators; braces
$\{\dots\}$ are always for Wronskians: $\{ f, g \}_D = gDf - fDg$,
where $D$ is a differential operator.

\subsection*{Acknowledgement}
The author is indebted to his advisor Professor Mark Adler, who
suggested me the problem and gave me technical helps. I also thank
him for valuable advice on writing and warm encouragement.

\section{The joint probability in the Dyson Brownian motion}
\label{PDE_for_the_finite_Dyson_process}

To get the PDE, we need to consider a generalized integral such that
indices $i$ and $j$ can be any positive integers
\begin{multline} \label{eq:definition_of_general_tau_n}
\tau_{n}(t^{(l)}_i, c^{(l)}_{i,j}; a^{(l)}_i) = \frac{1}{C}
\idotsint_{{U^{(0)}}^n \times \dots \times {U^{(m)}}^n}
V(\lambda^{(0)})
V(\lambda^{(m)}) \\
\prod^{m}_{l=0}
e^{\sum^{\infty}_{i=1}t^{(l)}_{i}\sum^{n}_{k=1}{\lambda^{(l)}_{k}}^{i}}
\prod^{m}_{l=1} e^{\sum^{\infty}_{i,j=1}
c^{(l)}_{i,j}\sum^{n}_{k=1}{\lambda^{(l-1)}_{k}}^{i}{\lambda^{(l)}_{k}}^{j}}
\prod^{m}_{l=0}\prod^{n}_{k=1}d\lambda^{(l)}_{k},
\end{multline}
with $C$ a normalization constant such that $\left. \PPDyson =
\tau_n \right|_{\mathfrak{L}}$, where the locus $\mathfrak{L}$ is
defined as ($l = 1, 2, \dots, m-1$, $k = 1, 2, \dots, m$, $c_k =
e^{-s_k}$)
\begin{equation} \label{eq:definition_of_locus}
\mathfrak{L} = \left\{
\begin{split}
t^{(0)}_{2} & = -\frac{1}{1-c^{2}_{1}}, \\
t^{(l)}_{2} & =
-\left( \frac{1}{1-c^2_l} + \frac{c^2_{l+1}}{1-c^2_{l+1}} \right), \\
t^{(m)}_{2} & = -\frac{1}{1-c^{2}_{m}}, \\
c^{(k)}_{1,1} & = \frac{2c_k}{1-c^2_k},  \\
& \textrm{and all other coefficients $0$.}
\end{split} \right.
\end{equation}
\begin{remark}
All $t^{(l)}_i$ and $c^{(l)}_{i,j}$ are variables of $\PPDyson$ in
latter part of the paper, though most of them assume the value $0$.
Therefore it is legitimate to consider $\frac{\partial}{\partial
t^{(0)}_1} \PPDyson$ etc.
\end{remark}

\begin{remark}
Since we allow $s^{(l)}_1$ to be $-\infty$ and $a^{(l)}_{2r_l}$ to
be $+\infty$, the integral in (\ref{eq:definition_of_general_tau_n})
may be divergent for general values of $t^{(l)}_{i}$ and
$c^{(k)}_{i,j}$. However, if we assume $t^{(l)}_{i} = 0$ for $i>2$ ,
$c^{(k)}_{i,j} = 0$ for $\max(i,j)>1$, and values of $t^{(l)}_{1}$,
$t^{(l)}_{2}$ and $c^{(k)}_{1,1}$ are near to the locus
$\mathfrak{L}$, then the integral is convergent, and all algebraic
operations in latter part of the paper can be taken in this
restricted setting, so they are legitimate.
\end{remark}

Now we consider actions of $D^{l,1}$ on $\tau_n$. Since $D^{l,1}$
acts on the integral domains of $\lambda^{(l)}_1, \dots,
\lambda^{(l)}_n$, by the formula
\begin{equation}
\left( \frac{\partial}{\partial a} + \frac{\partial}{\partial b}
\right) \int^b_a f(x)dx = f(b) - f(a) = \int^b_a f'(x)dx,
\end{equation}
we get
\begin{equation} \label{eq:action_of_D01}
\begin{split}
D^{0,1} \tau_n = & \frac{1}{C} \sum^{2r_0}_{i=1} \left[
\frac{\partial}{\partial a^{(0)}_i} \right] \idotsint_{{U^{(0)}}^n
\times
\dots \times {U^{(m)}}^n} V(\lambda^{(0)}) V(\lambda^{(m)}) \\
& \prod^{m}_{l=0}
e^{\sum^{\infty}_{i=1}t^{(l)}_{i}\sum^{n}_{k=1}{\lambda^{(l)}_{k}}^{i}}
\prod^{m}_{l=1} e^{\sum^{\infty}_{i,j=1}
c^{(l)}_{i,j}\sum^{n}_{k=1}{\lambda^{(l-1)}_{k}}^{i}{\lambda^{(l)}_{k}}^{j}}
\prod^{m}_{l=0}\prod^{n}_{k=1}d\lambda^{(l)}_{k} \\
= & \frac{1}{C} \idotsint_{{U^{(0)}}^n \times \dots \times
{U^{(m)}}^n} \left[ \sum^n_{k=1}\frac{\partial}{\partial
\lambda^{(0)}_k} \right] \left( V(\lambda^{(0)}) V(\lambda^{(m)})
\vphantom{\prod^{m}_{l=0}
e^{\sum^{\infty}_{i=1}t^{(l)}_{i}\sum^{n}_{k=1}{\lambda^{(l)}_{k}}^{i}}}
\right. \\
& \left. \prod^{m}_{l=0}
e^{\sum^{\infty}_{i=1}t^{(l)}_{i}\sum^{n}_{k=1}{\lambda^{(l)}_{k}}^{i}}
\prod^{m}_{l=1} e^{\sum^{\infty}_{i,j=1}
c^{(l)}_{i,j}\sum^{n}_{k=1}{\lambda^{(l-1)}_{k}}^{i}{\lambda^{(l)}_{k}}^{j}}
\right) \prod^{m}_{l=0}\prod^{n}_{k=1}d\lambda^{(l)}_{k} \\
= & \frac{1}{C} \idotsint_{{U^{(0)}}^n \times \dots \times
{U^{(m)}}^n} \left[ \sum^{\infty}_{i=1}it^{(0)}_i
\sum^n_{k=1}{\lambda^{(0)}_k}^{i-1} +
\sum^{\infty}_{i,j=0}ic^{(1)}_{i,j}
\sum^n_{k=1}{\lambda^{(0)}_k}^{i-1}{\lambda^{(1)}_k}^j \right] \\
& \left( V(\lambda^{(0)}) V(\lambda^{(m)}) \prod^{m}_{l=0}
e^{\sum^{\infty}_{i=1}t^{(l)}_{i}\sum^{n}_{k=1}{\lambda^{(l)}_{k}}^{i}}
\prod^{m}_{l=1} e^{\sum^{\infty}_{i,j=1}
c^{(l)}_{i,j}\sum^{n}_{k=1}{\lambda^{(l-1)}_{k}}^{i}{\lambda^{(l)}_{k}}^{j}}
\right) \prod^{m}_{l=0}\prod^{n}_{k=1}d\lambda^{(l)}_{k} \\
= & \frac{1}{C} \idotsint_{{U^{(0)}}^n \times \dots \times
{U^{(m)}}^n} \left[ nt^{(0)}_{1} + \sum^{\infty}_{i=2}
it^{(0)}_{i}\frac{\partial}{\partial t^{(0)}_{i-1}} +
\sum^{\infty}_{i=1} c^{(1)}_{1,i} \frac{\partial}{\partial
t^{(1)}_{i}} + \sum^{\infty}_{i=2} \sum^{\infty}_{j=1}
ic^{(1)}_{i,j}\frac{\partial}{\partial c^{(1)}_{i-1,j}}\right] \\
& \left( V(\lambda^{(0)}) V(\lambda^{(m)}) \prod^{m}_{l=0}
e^{\sum^{\infty}_{i=1}t^{(l)}_{i}\sum^{n}_{k=1}{\lambda^{(l)}_{k}}^{i}}
\prod^{m}_{l=1} e^{\sum^{\infty}_{i,j=1}
c^{(l)}_{i,j}\sum^{n}_{k=1}{\lambda^{(l-1)}_{k}}^{i}{\lambda^{(l)}_{k}}^{j}}
\right) \prod^{m}_{l=0}\prod^{n}_{k=1}d\lambda^{(l)}_{k} \\
= & \left[ nt^{(0)}_{1} + \sum^{\infty}_{i=2}
it^{(0)}_{i}\frac{\partial}{\partial t^{(0)}_{i-1}} +
\sum^{\infty}_{i=1} c^{(1)}_{1,i} \frac{\partial}{\partial
t^{(1)}_{i}} + \sum^{\infty}_{i=2} \sum^{\infty}_{j=1}
ic^{(1)}_{i,j}\frac{\partial}{\partial c^{(1)}_{i-1,j}}\right]
\tau_{n},
\end{split}
\end{equation}
and similarly ($l = 1, \dots, m-1$)
\begin{align}
D^{l,1}\tau_{n} = & \left[ nt^{(l)}_{1} + \sum^{\infty}_{i=2}
it^{(l)}_{i}\frac{\partial}{\partial t^{(l)}_{i-1}} +
\sum^{\infty}_{i=1} c^{(l)}_{i,1} \frac{\partial}{\partial
t^{(l-1)}_{i}} + \sum^{\infty}_{j=2} \sum^{\infty}_{i=1}
jc^{(l)}_{i,j}\frac{\partial}{\partial c^{(l)}_{i,j-1}}\right. \notag \\
& + \left. \sum^{\infty}_{i=1} c^{(l+1)}_{1,i}
\frac{\partial}{\partial t^{(l+1)}_{i}} + \sum^{\infty}_{i=2}
\sum^{\infty}_{j=1}
ic^{(l+1)}_{i,j}\frac{\partial}{\partial c^{(l+1)}_{i-1,j}}\right] \tau_{n}, \\
 D^{m,1}\tau_{n} = & \left[ nt^{(m)}_{1} +
\sum^{\infty}_{i=2} it^{(m)}_{i}\frac{\partial}{\partial
t^{(m)}_{i-1}} + \sum^{\infty}_{i=1} c^{(m)}_{i,1}
\frac{\partial}{\partial t^{(m-1)}_{i}} + \sum^{\infty}_{j=2}
\sum^{\infty}_{i=1} jc^{(m)}_{i,j}\frac{\partial}{\partial
c^{(m)}_{i,j-1}}\right] \tau_{n}. \label{eq:action_of_Dm1}
\end{align}
On the locus $\mathfrak{L}$ we get (We abbreviate $\PPDyson$ as
$\PPn$ here and latter)
\begin{align}
D^{0,1} \PPn = &
\left[-\frac{2}{1-c^{2}_{1}}\frac{\partial}{\partial t^{(0)}_{1}} +
\frac{2c_{1}}{1-c^{2}_{1}} \frac{\partial}{\partial
t^{(1)}_{1}}\right] \PPn, \label{eq:action_of_D01_on_locus} \\
D^{l,1} \PPn = & \left[ \frac{2c_{l}}{1-c^{2}_{l}}
\frac{\partial}{\partial t^{(l-1)}_{1}} -\left( \frac{2}{1-c^2_l} +
\frac{2c^2_{l+1}}{1-c^2_{l+1}} \right) \frac{\partial}{\partial
t^{(l)}_{1}} + \frac{2c_{l+1}}{1-c^{2}_{l+1}}
\frac{\partial}{\partial
t^{(l+1)}_{1}}\right] \PPn, \\
D^{m,1} \PPn = & \left[ \frac{2c_{m}}{1-c^{2}_{m}}
\frac{\partial}{\partial t^{(m-1)}_{1}} -
\frac{2}{1-c^{2}_{m}}\frac{\partial} {\partial t^{(m)}_{1}} \right]
\PPn. \label{eq:action_of_Dm1_on_locus}
\end{align}

Now we define an $m+1 \times m+1$ matrix
\begin{equation}
J = \begin{pmatrix}
2t^{(0)}_{2} & c^{(1)}_{1,1} & & & \\
c^{(1)}_{1,1} & 2t^{(1)}_{2} & c^{(2)}_{1,1} & & \\
& c^{(2)}_{1,1} & \ddots & \ddots & \\
& & \ddots & \ddots & c^{(m)}_{1,1} \\
& & & c^{(m)}_{1,1} & 2t^{(m)}_{2}
\end{pmatrix}^{-1},
\end{equation}
whose rows and columns are indexed from $0$ to $m$. On the locus
\begin{equation}
\left. J \right|_{\mathfrak{L}} = \begin{pmatrix}
-\frac{2}{1-c^{2}_{1}} & \frac{2c_{1}}{1-c^{2}_{1}} & & & \\
\frac{2c_{1}}{1-c^{2}_{1}} & -\frac{2}{1-c^2_l} -
\frac{2c^2_{l+1}}{1-c^2_{l+1}} &
\frac{2c_{2}}{1-c^{2}_{2}} & & \\
& \frac{2c_{2}}{1-c^{2}_{2}} & \ddots & \ddots & \\
& & \ddots & \ddots & \frac{2c_{m}}{1-c^{2}_{m}} \\
& & & \frac{2c_{m}}{1-c^{2}_{m}} & -\frac{2}{1-c^{2}_{m}} \\
\end{pmatrix}^{-1}.
\end{equation}
We can find the entries of the first and the last row of $J$ on the
locus explicitly,
\begin{align}
J_{0,l}|_{\mathfrak{L}} = & -\frac{1}{2}\prod^{l}_{i=1}c_{i} =
-\frac{1}{2}e^{-t_l}, \label{eq:first_row_of_J} \\
J_{m,l}|_{\mathfrak{L}} = & -\frac{1}{2}\prod^{m-l}_{i=1}c_{m-i+1} =
-\frac{1}{2}e^{t_l - t_m}, \label{eq:last_row_of_J} \\
\intertext{and especially} J_{0,m}|_{\mathfrak{L}} =
J_{m,0}|_{\mathfrak{L}} = & -\frac{1}{2}\prod^{m}_{i=1}c_{i} =
-\frac{1}{2}e^{-t_m}.
\end{align}

Then let
\begin{equation} \label{eq:definition_of_El1}
\begin{pmatrix}
E^{0,1} \\
E^{1,1} \\
\vdots \\
E^{m,1}
\end{pmatrix} = J
\begin{pmatrix}
D^{0,1} \\
D^{1,1} \\
\vdots \\
D^{m,1}
\end{pmatrix},
\end{equation}
we have
\begin{lemma} \label{lemma:e0emtau}
\begin{equation} \label{eq:e0emtau}
E^{0,1}E^{m,1}\log\PPn = E^{m,1}E^{0,1}\log\PPn = \frac{\partial^{2}
\log\PPn}{\partial t^{(0)}_{1} \partial t^{(m)}_{1}} -
\frac{n}{2}e^{-t_m}.
\end{equation}
\end{lemma}

\begin{proof}
First, since $E^{0,1}$ and $E^{m,1}$ are linear combinations of
$D^{l,1}$'s, they are differential operators of order $1$, and we
have
\begin{equation} \label{eq:log_identity}
E^{0,1}E^{m,1}\log\PPn = -\frac{E^{0,1}\PPn E^{m,1}\PPn}{\PPn^2} +
\frac{E^{0,1}E^{m,1}\PPn}{\PPn}.
\end{equation}
By (\ref{eq:action_of_D01_on_locus}) --
(\ref{eq:action_of_Dm1_on_locus}), (\ref{eq:first_row_of_J}) and
(\ref{eq:last_row_of_J}), we get
\begin{align}
E^{0,1}\PPn = & \left( \sum^m_{i=0}J_{0,i}|_{\mathfrak{L}}D^{i,1}
\right) \PPn =
\frac{\partial \PPn}{\partial t^{(0)}_1}, \\
\intertext{and} E^{m,1}\PPn = & \left( \sum^m_{i=0}J_{m,i}D^{i,1}
\right) \PPn = \frac{\partial \PPn}{\partial t^{(m)}_1}.
\end{align}
Therefore
\begin{equation} \label{eq:substitution_for_El1}
E^{0,1}E^{m,1}\log\PPn = -\frac{\frac{\partial \PPn}{\partial
t^{(0)}_1} \frac{\partial \PPn}{\partial t^{(m)}_1}}{\PPn^2} +
\frac{E^{0,1}\frac{\partial}{\partial t^{(m)}_1} \PPn}{\PPn}.
\end{equation}
Here we need to be careful about the term
$E^{0,1}\frac{\partial}{\partial t^{(m)}_1} \PPn$. By
(\ref{eq:action_of_D01}) -- (\ref{eq:action_of_Dm1}), the action of
$E^{m,1}$ on $\tau_n$ is equivalent to that of a differential
operator which does not contain $a^{(l)}_i$ explicitly. On the locus
$\mathfrak{L}$, all terms of the differential operator except for
$\frac{\partial}{\partial t^{(m)}_1}$ vanish, so we can ignore them
and replace $E^{m,1}$ by $\frac{\partial}{\partial t^{(m)}_1}$
between $E^{0,1}$ and $\PPn$.

Since $E^{0,1}$ and $\frac{\partial}{\partial t^{(m)}_1}$ commute,
\begin{equation} \label{eq:commutativity}
E^{0,1}\frac{\partial}{\partial t^{(m)}_1} \PPn =
\frac{\partial}{\partial t^{(m)}_1}E^{0,1} \PPn,
\end{equation}
by (\ref{eq:action_of_D01}) -- (\ref{eq:action_of_Dm1}),
(\ref{eq:first_row_of_J}) and (\ref{eq:definition_of_El1}), we have
the identity for the action of $E^{0,1}$ on $\tau_n$
\begin{equation} \label{eq:careful_substitution}
\begin{split}
E^{0,1}\tau_n = & J_{0,0} \left[ nt^{(0)}_{1} + \sum^{\infty}_{i=2}
it^{(0)}_{i}\frac{\partial}{\partial t^{(0)}_{i-1}} +
\sum^{\infty}_{i=1} c^{(1)}_{1,i} \frac{\partial}{\partial
t^{(1)}_{i}} + \sum^{\infty}_{i=2} \sum^{\infty}_{j=1}
ic^{(1)}_{i,j}\frac{\partial}{\partial c^{(1)}_{i-1,j}}\right]
\tau_{n} \\
& + \sum^{m-1}_{l=1} J_{0,l} \left[ nt^{(l)}_{1} +
\sum^{\infty}_{i=2} it^{(l)}_{i}\frac{\partial}{\partial
t^{(l)}_{i-1}} + \sum^{\infty}_{i=1} c^{(l)}_{i,1}
\frac{\partial}{\partial t^{(l-1)}_{i}} + \sum^{\infty}_{j=2}
\sum^{\infty}_{i=1}
jc^{(l)}_{i,j}\frac{\partial}{\partial c^{(l)}_{i,j-1}}\right. \\
& + \left. \sum^{\infty}_{i=1} c^{(l+1)}_{1,i}
\frac{\partial}{\partial t^{(l+1)}_{i}} + \sum^{\infty}_{i=2}
\sum^{\infty}_{j=1} ic^{(l+1)}_{i,j}\frac{\partial}{\partial
c^{(l+1)}_{i-1,j}}\right] \tau_{n} \\
& + J_{0,m} e^{-t_m} \left[ nt^{(m)}_{1} + \sum^{\infty}_{i=2}
it^{(m)}_{i}\frac{\partial}{\partial t^{(m)}_{i-1}} +
\sum^{\infty}_{i=1} c^{(m)}_{i,1} \frac{\partial}{\partial
t^{(m-1)}_{i}} + \sum^{\infty}_{j=2} \sum^{\infty}_{i=1}
jc^{(m)}_{i,j}\frac{\partial}{\partial
c^{(m)}_{i,j-1}}\right] \tau_{n} \\
= & \left( \frac{\partial}{\partial t^{(m)}_1} + nJ_{0,m}t^{(m)}_1 +
\dots \right) \tau_n,
\end{split}
\end{equation}
with coefficients of all terms except for $\frac{\partial}{\partial
t^{(m)}_1} + nJ_{0,m}t^{(m)}_1$ of the right-hand side operator
vanishing on the locus $\mathfrak{L}$ and not containing $t^{(m)}_1$
explicitly. So on the locus
\begin{equation} \label{eq:action_of_partial_tm1_E01}
\frac{\partial}{\partial t^{(m)}_1}E^{0,1} \PPn =
\frac{\partial}{\partial t^{(m)}_1} \left( \frac{\partial}{\partial
t^{(m)}_1} + nJ_{0,m}t^{(m)}_1 + \dots \right) \PPn =
\frac{\partial^2 \PPn}{\partial t^{(0)}_1
\partial t^{(m)}_1} + nJ_{0,m}|_{\mathfrak{L}} \PPn,
\end{equation}
and
\begin{equation}
E^{0,1}E^{m,1}\log\PPn = -\frac{\frac{\partial \PPn}{\partial
t^{(0)}_1} \frac{\partial \PPn}{\partial t^{(m)}_1}}{\PPn^2} +
\frac{\frac{\partial^2 \PPn}{\partial t^{(0)}_1
\partial t^{(m)}_1} + nJ_{0,m}|_{\mathfrak{L}} \PPn}{\PPn} = \frac{\partial^{2}
\log\PPn}{\partial t^{(0)}_{1} \partial t^{(m)}_{1}} -
\frac{n}{2}e^{-t_m}.
\end{equation}
\end{proof}

Similarly to (\ref{eq:action_of_D01}) -- (\ref{eq:action_of_Dm1}),
with the help of the formula
\begin{equation}
\left( a\frac{\partial}{\partial a} + b\frac{\partial}{\partial b}
\right) \int^b_a f(x)dx = bf(b) - af(a) = \int^b_a (xf(x))'dx,
\end{equation}
we get ($\left[ \sum^n_{k=1} \frac{\partial}{\partial
\lambda^{(0)}_k} \lambda^{(0)}_k \right]$ is regarded as an
operator)
\begin{equation} \label{eq:action_of_D02}
\begin{split}
D^{0,2} \tau_n = & \frac{1}{C} \sum^{2r_0}_{i=1} \left[ a^{(0)}_i
\frac{\partial}{\partial a^{(0)}_i} \right] \idotsint_{{U^{(0)}}^n
\times
\dots \times {U^{(m)}}^n} V(\lambda^{(0)}) V(\lambda^{(m)}) \\
& \prod^{m}_{l=0}
e^{\sum^{\infty}_{i=1}t^{(l)}_{i}\sum^{n}_{k=1}{\lambda^{(l)}_{k}}^{i}}
\prod^{m}_{l=1} e^{\sum^{\infty}_{i,j=1}
c^{(l)}_{i,j}\sum^{n}_{k=1}{\lambda^{(l-1)}_{k}}^{i}{\lambda^{(l)}_{k}}^{j}}
\prod^{m}_{l=0}\prod^{n}_{k=1}d\lambda^{(l)}_{k} \\
= & \frac{1}{C} \idotsint_{{U^{(0)}}^n \times \dots \times
{U^{(m)}}^n} \left[ \sum^n_{k=1}\frac{\partial}{\partial
\lambda^{(0)}_k} \lambda^{(0)}_k \right] \left( V(\lambda^{(0)})
V(\lambda^{(m)})
\vphantom{\prod^{m}_{l=0}
e^{\sum^{\infty}_{i=1}t^{(l)}_{i}\sum^{n}_{k=1}{\lambda^{(l)}_{k}}^{i}}}
\right. \\
& \left. \prod^{m}_{l=0}
e^{\sum^{\infty}_{i=1}t^{(l)}_{i}\sum^{n}_{k=1}{\lambda^{(l)}_{k}}^{i}}
\prod^{m}_{l=1} e^{\sum^{\infty}_{i,j=1}
c^{(l)}_{i,j}\sum^{n}_{k=1}{\lambda^{(l-1)}_{k}}^{i}{\lambda^{(l)}_{k}}^{j}}
\right) \prod^{m}_{l=0}\prod^{n}_{k=1}d\lambda^{(l)}_{k} \\
= & \frac{1}{C} \idotsint_{{U^{(0)}}^n \times \dots \times
{U^{(m)}}^n} \left[ \sum^{\infty}_{i=1}it^{(0)}_i
\sum^n_{k=1}{\lambda^{(0)}_k}^i +
\sum^{\infty}_{i,j=0}ic^{(1)}_{i,j}
\sum^n_{k=1}{\lambda^{(0)}_k}^i{\lambda^{(1)}_k}^j + \frac{n(n+1)}{2} \right] \\
& \left( V(\lambda^{(0)}) V(\lambda^{(m)}) \prod^{m}_{l=0}
e^{\sum^{\infty}_{i=1}t^{(l)}_{i}\sum^{n}_{k=1}{\lambda^{(l)}_{k}}^{i}}
\prod^{m}_{l=1} e^{\sum^{\infty}_{i,j=1}
c^{(l)}_{i,j}\sum^{n}_{k=1}{\lambda^{(l-1)}_{k}}^{i}{\lambda^{(l)}_{k}}^{j}}
\right) \prod^{m}_{l=0}\prod^{n}_{k=1}d\lambda^{(l)}_{k} \\
= & \frac{1}{C} \idotsint_{{U^{(0)}}^n \times \dots \times
{U^{(m)}}^n} \left[ \sum^{\infty}_{i=2}
it^{(0)}_{i}\frac{\partial}{\partial t^{(0)}_i} +
\sum^{\infty}_{i,j=1}
ic^{(1)}_{i,j}\frac{\partial}{\partial c^{(1)}_{i,j}} + \frac{n(n+1)}{2} \right] \\
& \left( V(\lambda^{(0)}) V(\lambda^{(m)}) \prod^{m}_{l=0}
e^{\sum^{\infty}_{i=1}t^{(l)}_{i}\sum^{n}_{k=1}{\lambda^{(l)}_{k}}^{i}}
\prod^{m}_{l=1} e^{\sum^{\infty}_{i,j=1}
c^{(l)}_{i,j}\sum^{n}_{k=1}{\lambda^{(l-1)}_{k}}^{i}{\lambda^{(l)}_{k}}^{j}}
\right) \prod^{m}_{l=0}\prod^{n}_{k=1}d\lambda^{(l)}_{k} \\
= & \left[\sum^{\infty}_{i=1} it^{(0)}_{i}\frac{\partial}{\partial
t^{(0)}_{i}} + \sum^{\infty}_{i=1} \sum^{\infty}_{j=1}
ic^{(1)}_{i,j} \frac{\partial}{\partial c^{(1)}_{i,j}} +
\frac{n(n+1)}{2}\right]\tau_{n},
\end{split}
\end{equation}
and similarly
\begin{align}
D^{m,2}\tau_{n} = & \left[\sum^{\infty}_{i=1}
it^{(m)}_{i}\frac{\partial}{\partial t^{(m)}_{i}} +
\sum^{\infty}_{i=1} \sum^{\infty}_{j=1} jc^{(m)}_{i,j}
\frac{\partial}{\partial c^{(m)}_{i,j}} +
\frac{n(n+1)}{2}\right]\tau_{n}, \\
\intertext{and for $l = 1, \dots, m-1$} D^{l,2}\tau_{n} = &
\left[\sum^{\infty}_{i=1} it^{(l)}_{i}\frac{\partial}{\partial
t^{(l)}_{i}} + \sum^{\infty}_{i=1} \sum^{\infty}_{j=1}
jc^{(l)}_{i,j} \frac{\partial}{\partial c^{(l)}_{i,j}} +
\sum^{\infty}_{i=1} \sum^{\infty}_{j=1} ic^{(l+1)}_{i,j}
\frac{\partial}{\partial c^{(l+1)}_{i,j}} + n \right]\tau_{n}.
\end{align}
On the locus $\mathfrak{L}$ we get ($l = 1, \dots, m-1$)
\begin{align}
D^{0,2}\PPn = & \left[-\frac{2}{1-c^{2}_{1}}
\frac{\partial}{\partial t^{(0)}_{2}} + \frac{2c_{1}}{1-c^{2}_{1}}
\frac{\partial}{\partial c^{(1)}_{1,1}} + \frac{n(n+1)}{2}\right] \PPn, \\
D^{l,2}\PPn = & \left[ -\left( \frac{2}{1-c^2_l} +
\frac{c^2_{l+1}}{1-c^2_{l+1}} \right) \frac{\partial}{\partial
t^{(l)}_{2}} + \frac{2c_{l}}{1-c^{2}_{l}} \frac{\partial}{\partial
c^{(l)}_{1,1}} + \frac{2c_{l+1}}{1-c^{2}_{l+1}}
\frac{\partial}{\partial
c^{(l+1)}_{1,1}} + n \right] \PPn, \\
D^{m,2}\PPn = & \left[-\frac{2}{1-c^{2}_{m}}
\frac{\partial}{\partial t^{(m)}_{2}} + \frac{2c_{m}}{1-c^{2}_{m}}
\frac{\partial}{\partial c^{(m)}_{1,1}} + \frac{n(n+1)}{2}\right]
\PPn.
\end{align}

If we define ($l = 1, 2, \dots, m-1$)
\begin{align}
E^{0,2} = & D^{0,2} - c^{(1)}_{1,1} \frac{\partial}{\partial
c^{(1)}_{1,1}} , \label{eq:definition_of_E02} \\
E^{l,2} = & D^{l,2} - c^{(l)}_{1,1} \frac{\partial}{\partial
c^{(l)}_{1,1}} - - c^{(l+1)}_{1,1} \frac{\partial}{\partial
c^{(l+1)}_{1,1}} , \\
E^{m,2} = & D^{m,2} - c^{(m)}_{1,1} \frac{\partial}{\partial
c^{(m)}_{1,1}} , \label{eq:definition_of_Em2}
\end{align}
we have
\begin{lemma}
For $k, l = 0, 1, \dots, m$,
\begin{equation} \label{eq:ekeltau}
E^{k,2}E^{l,1}\log\PPn = 2t^{(k)}_{2}|_{\mathfrak{L}}
\frac{\partial^{2} \PPn}{\partial t^{(k)}_{2} \partial t^{(l)}_{1}}
+ \delta^{l}_{k}E^{l,1} \log\PPn.
\end{equation}
\end{lemma}

\begin{proof}
With arguments similar to those for (\ref{eq:log_identity}) and
(\ref{eq:substitution_for_El1}), we get
\begin{equation}
E^{k,2}E^{l,1}\log\PPn = \frac{E^{k,2}\PPn E^{l,1}\PPn}{\PPn^2} +
\frac{E^{k,2}E^{l,1}\PPn}{\PPn} = -\frac{\left(
2t^{(k)}_2|_{\mathfrak{L}} \frac{\partial \PPn}{\partial t^{(k)}_2}
+ C\PPn \right) \frac{\partial \PPn}{\partial t^{(l)}_1}}{\PPn^2} +
\frac{E^{k,2} \frac{\partial}{\partial t^{(l)}_1} \PPn}{\PPn}.
\end{equation}
Here
\begin{equation}
C=
\begin{cases}
-\frac{n(n+1)}{2} & k = 1 \text{ or } m, \\
-n & \text{otherwise.}
\end{cases}
\end{equation}

Similar to (\ref{eq:commutativity}) and
(\ref{eq:careful_substitution}), we have
\begin{equation}
E^{k,2} \frac{\partial}{\partial t^{(l)}_1} \PPn =
\frac{\partial}{\partial t^{(l)}_1}E^{k,2} \PPn
\end{equation}
and
\begin{equation}
E^{k,2}\tau_n = \left( 2t^{(k)}_2\frac{\partial}{\partial t^{(k)}_2}
+ t^{(k)}_1 \frac{\partial}{\partial t^{(k)}_1} + C + \dots \right)
\tau_n,
\end{equation}
with coefficients of all terms other than
$2t^{(k)}_2\frac{\partial}{\partial t^{(k)}_2}$, $t^{(k)}_1
\frac{\partial}{\partial t^{(k)}_1}$ or $C$ of the right-hand side
operator vanishing on $\mathfrak{L}$ and not containing $t^{(l)}_1$
explicitly. Therefore with an argument similar to that for
(\ref{eq:action_of_partial_tm1_E01})
\begin{equation}
\begin{split}
E^{k,2}E^{l,1}\log\PPn = & -\frac{\left( 2t^{(k)}_2|_{\mathfrak{L}}
\frac{\partial \PPn}{\partial t^{(k)}_2} + C\PPn \right)
\frac{\partial \PPn}{\partial t^{(l)}_1}}{\PPn^2} +
\frac{\frac{\partial}{\partial t^{(l)}_1} \left(
2t^{(k)}_2\frac{\partial}{\partial t^{(k)}_2} + t^{(k)}_1
\frac{\partial}{\partial t^{(k)}_1} + C \right)
\PPn}{\PPn} \\
= & 2t^{(k)}_2|_{\mathfrak{L}} \left( -\frac{\frac{\partial
\PPn}{\partial t^{(k)}_2} \frac{\partial \PPn}{\partial
t^{(l)}_1}}{\PPn^2} + \frac{\frac{\partial^2 \PPn}{\partial
t^{(k)}_2
\partial^{(l)}_1}}{\PPn} \right) + \frac{\delta^l_k \frac{\partial \PPn}{\partial
t^{(l)}_1} + t^{(k)}_1|_{\mathfrak{L}} \frac{\partial^2 \PPn}
{\partial t^{(k)}_1 \partial t^{(l)}_1}}{\PPn} \\
= & 2t^{(k)}_2|_{\mathfrak{L}} \frac{\partial^2 \log\PPn}{\partial
t^{(k)}_2
\partial^{(l)}_1} + \delta^l_k E^{l,1} \log\PPn,
\end{split}
\end{equation}
since
\begin{equation}
E^{l,1} \log\PPn = \frac{E^{l,1} \PPn}{\PPn} = \frac{\frac{\partial
\PPn}{\partial t^{(l)}_1}}{\PPn}.
\end{equation}
\end{proof}

Since on the locus $\mathfrak{L}$, $c^{(k)}_{1,1}$ and $t^{(l)}_{2}$
are functions of $c_1 = e^{-s_1}, \dots, c_{m} = e^{-s_m}$ defined
in (\ref{eq:definition_of_locus}), by the chain rule we get as
operators on $\PPn$ ($l = 1, 2, \dots, m$)
\begin{equation}
\frac{\partial}{\partial s_{l}} = \frac{2c^2_l}{(1-c^2_l)^{2}}
\frac{\partial}{\partial t^{(l-1)}_{2}} +
\frac{2c^2_l}{(1-c^2_l)^{2}} \frac{\partial}{\partial t^{(l)}_{2}} -
\frac{2c_l(1+c^2_l)}{(1-c^2_l)^{2}} \frac{\partial}{\partial
c^{(l)}_{1,1}}
\end{equation}
and
\begin{equation}
c^{(l)}_{1,1} \frac{\partial}{\partial c^{(l)}_{1,1}} =
\frac{2c^2_l}{1-c^4_l} \frac{\partial}{\partial t^{(l-1)}_{2}} +
\frac{2c^2_l}{1-c^4_l} \frac{\partial}{\partial t^{(l)}_{2}} -
\frac{1-c^2_l}{1+c^2_l} \frac{\partial}{\partial s_{l}}.
\end{equation}
Therefore by (\ref{eq:definition_of_E02}) --
(\ref{eq:definition_of_Em2}) we get on $\mathfrak{L}$ that ($l = 1,
2, \dots, m-1$)
\begin{align}
E^{0,2} = & D^{0,2} - \frac{2c^2_1}{1-c^4_1}
\frac{\partial}{\partial t^{(0)}_{2}} - \frac{2c^2_1}{1-c^4_1}
\frac{\partial}{\partial t^{(1)}_{2}} + \frac{1-c^2_1}{1+c^2_1}
\frac{\partial}{\partial s_{1}}, \label{eq:another_form_of_E02} \\
E^{l,2} = & D^{l,2} - \frac{2c^2_l}{1-c^4_l}
\frac{\partial}{\partial t^{(l-1)}_{2}} -
\left(\frac{2c^2_l}{1-c^4_l} +
\frac{2c^2_{l+1}}{1-c^4_{l+1}}\right)\frac{\partial}{\partial
t^{(l)}_{2}} \notag \\
& - \frac{2c^2_{l+1}}{1-c^4_{l+1}} \frac{\partial}{\partial
t^{(l+1)}_{2}} + \frac{1-c^2_l}{1+c^2_l} \frac{\partial}{\partial
s_{l}} + \frac{1-c^2_{l+1}}{1+c^2_{l+1}} \frac{\partial}{\partial
s_{l+1}}, \\
E^{m,2} = & D^{m,2} - \frac{2c^2_m}{1-c^4_m}
\frac{\partial}{\partial t^{(m-1)}_{2}} - \frac{2c^2_m}{1-c^4_m}
\frac{\partial}{\partial t^{(m)}_{2}} + \frac{1-c^2_m}{1+c^2_m}
\frac{\partial}{\partial s_{m}}.
\end{align}
Now we denote ($l = 1, 2, \dots, m-1$)
\begin{align}
F^{0,2} = & D^{0,2} + \frac{1-c^2_1}{1+c^2_1}
\frac{\partial}{\partial s_{1}}, \\
F^{l,2} = & D^{l,2}  + \frac{1-c^2_l}{1+c^2_l}
\frac{\partial}{\partial s_{l}} + \frac{1-c^2_{l+1}}{1+c^2_{l+1}}
\frac{\partial}{\partial s_{l+1}}, \\
F^{m,2} = & D^{m,2} + \frac{1-c^2_m}{1+c^2_m}
\frac{\partial}{\partial s_{m}},
\end{align}
and we have
\begin{lemma}
For $l = 1, \dots, m-1$,
\begin{align}
F^{0,2}E^{m,1} \log \PPn = &
\left[-\frac{2}{1-c^4_1}\frac{\partial^{2}}{\partial t^{(0)}_{2}
\partial t^{(m)}_{1}} + \frac{2c^2_1}{1-c^4_1}
\frac{\partial^{2}}{\partial
t^{(1)}_{2} \partial t^{(m)}_{1}}\right] \log \PPn, \label{eq:formula_of_F02Em1} \\
F^{l,2}E^{m,1} \log \PPn = &
\left[\frac{2c^2_l}{1-c^4_l}\frac{\partial^{2}}{\partial
t^{(l-1)}_{2} \partial t^{(m)}_{1}} - \left( \frac{2}{1-c^4_l} +
\frac{2c^4_{l+1}}{1-c^4_{l+1}} \right) \frac{\partial^{2}}{\partial
t^{(l)}_{2} \partial t^{(m)}_{1}} \right. \notag \\
& \left. +
\frac{2c^2_{l+1}}{1-c^4_{l+1}}\frac{\partial^{2}}{\partial
t^{(l+1)}_{2} \partial t^{(m)}_{1}} \right] \log \PPn, \\
F^{m,2}E^{m,1} \log \PPn = & \left[\frac{2c^2_m}{1-c^4_m}
\frac{\partial^{2}}{\partial t^{(m-1)}_{2} \partial t^{(m)}_{1}}
-\frac{2}{1-c^4_m}\frac{\partial^{2}}{\partial t^{(m)}_{2}
\partial t^{(m)}_{1}} \right] \log \PPn + E^{m,1} \log \PPn.
\label{eq:formula_of_Fm2Em1}
\end{align}
\end{lemma}

\begin{proof}
We only prove (\ref{eq:formula_of_F02Em1}). By
(\ref{eq:another_form_of_E02}), on $\mathfrak{L}$ we have
\begin{equation}
F^{0,2} = E^{0,2} + \frac{2c^2_1}{1-c^4_1} \frac{\partial}{\partial
t^{(0)}_{2}} + \frac{2c^2_1}{1-c^4_1} \frac{\partial}{\partial
t^{(1)}_{2}},
\end{equation}
and similar to (\ref{eq:substitution_for_El1}) and
(\ref{eq:action_of_partial_tm1_E01}), we have
\begin{equation}
\frac{\partial}{\partial t^{(k)}_2} E^{l,1} \log\PPn =
-\frac{\frac{\partial \PPn}{\partial t^{(k)}_2} E^{l,1}\PPn}{\PPn^2}
+ \frac{\frac{\partial}{\partial t^{(k)}_2} E^{l,1} \PPn}{\PPn} =
\frac{\partial^2 \log\PPn}{\partial t^{(k)}_2 \partial t^{(l)}_1},
\end{equation}
so that with the result of (\ref{eq:ekeltau}),
\begin{equation}
\begin{split}
F^{0,2}E^{m,1} \log \PPn = & E^{0,2}E^{m,1}\log\PPn +
\frac{2c^2_1}{1-c^4_1} \frac{\partial}{\partial
t^{(0)}_{2}}E^{m,1}\log\PPn + \frac{2c^2_1}{1-c^4_1}
\frac{\partial}{\partial t^{(1)}_{2}}E^{m,1}\log\PPn \\
= & 2t^{(0)}_2 \frac{\partial^2 \log\PPn}{\partial t^{(0)}_2
\partial t^{(m)}_1} + \frac{2c^2_1}{1-c^4_1}
\frac{\partial^2 \log\PPn}{\partial t^{(0)}_2 \partial t^{(m)}_1} +
\frac{2c^2_1}{1-c^4_1} \frac{\partial^2
\log\PPn}{\partial t^{(1)}_2 \partial t^{(m)}_1} \\
= & \left[-\frac{2}{1-c^4_1}\frac{\partial^{2}}{\partial t^{(0)}_{2}
\partial t^{(m)}_{1}} + \frac{2c^2_1}{1-c^4_1}
\frac{\partial^{2}}{\partial t^{(1)}_{2} \partial
t^{(m)}_{1}}\right] \log \PPn.
\end{split}
\end{equation}
\end{proof}

Finally we define
\begin{equation}
\begin{pmatrix}
G^{0,2}\\
G^{1,2}\\
\vdots \\
G^{m,2}\\
\end{pmatrix} = K
\begin{pmatrix}
F^{0,2}\\
F^{1,2}\\
\vdots \\
F^{m,2}\\
\end{pmatrix},
\end{equation}
where
\begin{equation}
K =
\begin{pmatrix}
-\frac{2}{1-c^4_1} & \frac{2c^2_1}{1-c^4_1} & & &
\\
\frac{2c^2_1}{1-c^4_1} & - \frac{2}{1-c^4_1} -
\frac{2c^4_2}{1-c^4_2} &
\frac{2c^2_2}{1-c^4_2} & & \\
& \frac{2c^2_2}{1-c^4_2} & \ddots & \ddots & \\
& & \ddots & \ddots & \frac{2c^2_m}{1-c^4_m} \\
& & & \frac{2c^2_m}{1-c^4_m} & -\frac{2}{1-c^4_m} \\
\end{pmatrix}^{-1}.
\end{equation}
We can get $K^{-1}$ by substituting each $c_i$ in
$J^{-1}|_{\mathfrak{L}}$ by $c^2_i$, so we have
\begin{equation}
K_{0,l} = -\frac{1}{2}\prod^{l}_{i=1}c^{2}_{i}, \qquad K_{m,l} =
-\frac{1}{2}\prod^{m-l}_{i=1}c^{2}_{m-i+1}.
\end{equation}

and get by (\ref{eq:formula_of_F02Em1}) --
(\ref{eq:formula_of_Fm2Em1}),

\begin{lemma}
\begin{equation}
G^{0,2}E^{m,1} \log \PPn = \frac{\partial^{2}}{\partial t^{(0)}_{2}
\partial t^{(m)}_{1}} \log \PPn + K_{0,m}E^{m,1} \log \PPn.
\end{equation}
\end{lemma}
Symmetrically, we can get by the same method
\begin{lemma}
\begin{equation}
G^{m,2}E^{1,1} \log \PPn = \frac{\partial^{2}}{\partial t^{(m)}_{2}
\partial t^{(0)}_{1}} \log \PPn + K_{m,0}E^{0,1} \log \PPn.
\end{equation}
\end{lemma}

By the result of \cite{Adler-van_Moerbeke99},
\begin{align}
\frac{\partial}{\partial t^{(0)}_{1}} \log
\frac{\tau_{n+1}}{\tau_{n-1}} = & \frac{\frac{\partial^{2}}{\partial
t^{(0)}_{2} \partial t^{(m)}_{1}} \log \tau_{n}}
{\frac{\partial^{2}}{\partial t^{(0)}_{1}
\partial t^{(m)}_{1}} \log \tau_{n}}, \\
\frac{\partial}{\partial t^{(m)}_{1}} \log
\frac{\tau_{n+1}}{\tau_{n-1}} = & \frac{\frac{\partial^{2}}{\partial
t^{(m)}_{2} \partial t^{(0)}_{1}} \log \tau_{n}}
{\frac{\partial^{2}}{\partial t^{(m)}_{1}
\partial t^{(0)}_{1}} \log \tau_{n}},
\end{align}
we get the differential equation on the locus
\begin{align}
\left. E^{0,1} \log
\frac{\tau_{n+1}}{\tau_{n-1}}\right|_{\mathfrak{L}} = & \frac{\left.
G^{0,2}E^{m,1} \log \tau_{n} \right|_{\mathfrak{L}} - \left.
K_{0,m}E^{m,1} \log \tau_{n} \right|_{\mathfrak{L}}} {\left.
E^{0,1}E^{m,1} \log \tau_{n} \right|_{\mathfrak{L}} - nJ_{m,0}}, \\
\left. E^{m,1} \log
\frac{\tau_{n+1}}{\tau_{n-1}}\right|_{\mathfrak{L}} = & \frac{\left.
G^{m,2}E^{0,1} \log \tau_{n} \right|_{\mathfrak{L}} - \left.
K_{m,0}E^{0,1} \log \tau_{n} \right|_{\mathfrak{L}}} {\left.
E^{m,1}E^{0,1} \log \tau_{n} \right|_{\mathfrak{L}} - nJ_{0,m}},
\end{align}
where the result
\begin{equation}
\left. \frac{\partial}{\partial t^{(0)}_{1}} \log
\frac{\tau_{n+1}}{\tau_{n-1}} \right|_{\mathfrak{L}} = \left.
E^{0,1} \log \frac{\tau_{n+1}}{\tau_{n-1}}\right|_{\mathfrak{L}}
\end{equation}
can be proved like lemma \ref{lemma:e0emtau}. By the identity
\begin{equation}
\left. E^{0,1}E^{m,1} \log
\frac{\tau_{n+1}}{\tau_{n-1}}\right|_{\mathfrak{L}} = \left.
E^{m,1}E^{0,1} \log
\frac{\tau_{n+1}}{\tau_{n-1}}\right|_{\mathfrak{L}},
\end{equation}
we get the final result
\begin{equation} \label{eq:equation for Dyson process}
E^{0,1}\frac{G^{m,2}E^{0,1} \log \tilde{\tau}_{n} - K_{m,0}E^{0,1}
\log \tilde{\tau}_{n}} {E^{m,1}E^{0,1} \log \tilde{\tau}_{n} -
nJ_{0,m}} = E^{m,1}\frac{G^{0,2}E^{m,1} \log \tilde{\tau}_{n} -
K_{0,m}E^{m,1} \log \tilde{\tau}_{n}} {E^{0,1}E^{m,1} \log
\tilde{\tau}_{n} - nJ_{m,0}}.
\end{equation}
Now we denote
\begin{align}
\mathcal{A}_1 = & -2E^{0,1}, \\
\mathcal{B}_1 = & -2E^{m,1}, \\
\mathcal{A}_2 = & -2(G^{0,2}-K_{0,m}), \\
\mathcal{B}_2 = & -2(G^{m,2}-K_{m,0}),
\end{align}
and we get the equation (\ref{eq:compact equation for Dyson
process}).

\begin{remark}
In the $2$-time case, i.e., $m=1$,
\begin{align}
\mathcal{A}_1 = & D^{0,1} + c_1 D^{1,1}, \\
\mathcal{B}_1 = & c_1 D^{0,1} + D^{1,1},\\
\mathcal{A}_2 = & F^{0,1} + c^2_1 F^{1,1} = D^{0,1} + c^2_1 D^{1,1}
+ (1-c^2_1)\frac{\partial}{\partial
t_1} -c^2_1, \\
\mathcal{B}_2 = & c^2_1 F^{0,1} + F^{1,1} = c^2_1 D^{0,1} + D^{1,1}
+ (1-c^2_1)\frac{\partial}{\partial t_1} -c^2_1.
\end{align}
Our PDE (\ref{eq:compact equation for Dyson process}) agrees with
that in \cite{Adler-van_Moerbeke05}.
\end{remark}

\section{The joint probability in the Airy process} \label{PDE_for_the_Airy_process}

In this section we adapt notations defined in (\ref{eq:first scaled
variable for Airy process})--(\ref{eq:last scaled variable for Airy
process}), and by remark \ref{remark:on_al1}, $a^{(l)}_1 =
\bar{a}^{(l)}_1 = -\infty$. We denote ($l = 0, 1, \dots, m$)
\begin{equation}
\bar{D}^{l,1} = \sum^{2r_{l}}_{k=1}\frac{\partial}{\partial
\bar{a}^{(l)}_{k}}, \qquad \bar{D}^{l,2} =
\sum^{2r_{l}}_{k=1}\bar{a}^{(l)}_{k}\frac{\partial}{\partial
\bar{a}^{(l)}_{k}},
\end{equation}
if all $a^{(l)}_{r_l} < +\infty$, otherwise drop the $a^{(l)}_{r_l}$
part. we can write the differential operators defined for the Dyson
process as
\begin{align}
\mathcal{A}_1 = & \sqrt{2\bar{n}} \sum^{m}_{l=0}
e^{-\bar{t}_l/\bar{n}} \bar{D}^{l,1},
\label{eq:first scaled operator for Airy process} \\
\mathcal{B}_1 = & \sqrt{2\bar{n}} \sum^{m}_{l=0}
e^{(\bar{t}_l - \bar{t}_m)/\bar{n}} \bar{D}^{l,1}, \\
\mathcal{A}_2 = & \sum^{m}_{l=0} e^{-2\bar{t}_l/\bar{n}}
\bar{D}^{l,2} +2\bar{n}^{2} \sum^{m}_{l=0} e^{-2\bar{t}_l/\bar{n}}
\bar{D}^{l,1} + \bar{n} \sum^{m}_{l=1} (1 - e^{-2\bar{t}_l/\bar{n}})
\frac{\partial} {\partial \bar{t}_{l}} - e^{-2\bar{t}_m/\bar{n}}, \\
\mathcal{B}_2 = & \sum^{m}_{l=0} e^{2(\bar{t}_l -
\bar{t}_m)/\bar{n}} \bar{D}^{l,2} + 2\bar{n}^{2} \sum^{m}_{l=0}
e^{2(\bar{t}_l - \bar{t}_m)/\bar{n}} \bar{D}^{l,1} + \bar{n}
\sum^{m}_{l=1} (e^{2(\bar{t}_l - \bar{t}_m)/\bar{n}} -
e^{-2\bar{t}_m/\bar{n}}) \frac{\partial} {\partial \bar{t}_{l}} -
e^{-2\bar{t}_m/\bar{n}}. \label{eq:last scaled operator for Airy
process}
\end{align}

It is not difficult to see that (\ref{eq:compact equation for Dyson
process}) implies
\begin{multline} \label{eq:new_equation_for_Dyson_process}
[\mathcal{A}_1 \mathcal{B}_2 \mathcal{A}_1 - \mathcal{B}_1
\mathcal{A}_2 \mathcal{B}_1] \log\PPn \cdot (\mathcal{A}_1
\mathcal{B}_1 \log\PPn + 2ne^{-t_m}) = \\
\mathcal{B}_2 \mathcal{A}_1 \log\PPn \cdot \mathcal{A}_1
\mathcal{B}_1 \mathcal{A}_1 \log\PPn - \mathcal{A}_2 \mathcal{B}_1
\log\PPn \cdot \mathcal{B}_1 \mathcal{A}_1 \mathcal{B}_1 \log\PPn.
\end{multline}

Take substitutions (\ref{eq:first scaled operator for Airy
process})--(\ref{eq:last scaled operator for Airy process}) into
(\ref{eq:new_equation_for_Dyson_process}), we get
\begin{equation} \label{eq:scaled_equation_for_Dyson_process_before_asymptotics}
\begin{split}
& \left( \left[ \sum^{m}_{l=0} e^{-\bar{t}_l/\bar{n}} \bar{D}^{l,1}
\right] \left[\sum^{m}_{l=0} e^{2(\bar{t}_l - \bar{t}_m)/\bar{n}}
\bar{D}^{l,2} + 2\bar{n}^{2} \sum^{m}_{l=0} e^{2(\bar{t}_l -
\bar{t}_m)/\bar{n}} \bar{D}^{l,1} \right. \right. \\
& \left. + \bar{n} \sum^{m}_{l=1} (e^{2(\bar{t}_l -
\bar{t}_m)/\bar{n}} - e^{-2\bar{t}_m/\bar{n}}) \frac{\partial}
{\partial \bar{t}_{l}} - e^{-2\bar{t}_m/\bar{n}} \right] \left[
\sum^{m}_{l=0} e^{-\bar{t}_l/\bar{n}} \bar{D}^{l,1} \right]
\log\PPn \\
& - \left[ \sum^{m}_{l=0} e^{(\bar{t}_l - \bar{t}_m)/\bar{n}}
\bar{D}^{l,1} \right] \left[ \sum^{m}_{l=0} e^{-2\bar{t}_l/\bar{n}}
\bar{D}^{l,2} +2\bar{n}^{2} \sum^{m}_{l=0} e^{-2\bar{t}_l/\bar{n}}
\bar{D}^{l,1} \right. \\
& \left. + \left. \bar{n} \sum^{m}_{l=1} (1 -
e^{-2\bar{t}_l/\bar{n}}) \frac{\partial} {\partial \bar{t}_{l}} -
e^{-2\bar{t}_m/\bar{n}} \right] \left[ \sum^{m}_{l=0} e^{(\bar{t}_l
- \bar{t}_m)/\bar{n}} \bar{D}^{l,1} \right] \log\PPn \right) \\
& \times \left( \left[ \sum^{m}_{l=0} e^{-\bar{t}_l/\bar{n}}
\bar{D}^{l,1} \right] \left[ \sum^{m}_{l=0} e^{(\bar{t}_l -
\bar{t}_m)/\bar{n}} \bar{D}^{l,1} \right] \log\PPn + \bar{n}^2
e^{-\bar{t}_m/\bar{n}} \right) \\
= & \left[ \sum^{m}_{l=0} e^{2(\bar{t}_l - \bar{t}_m)/\bar{n}}
\bar{D}^{l,2} + 2\bar{n}^{2} \sum^{m}_{l=0} e^{2(\bar{t}_l -
\bar{t}_m)/\bar{n}} \bar{D}^{l,1} + \bar{n} \sum^{m}_{l=1}
(e^{2(\bar{t}_l - \bar{t}_m)/\bar{n}} - e^{-2\bar{t}_m/\bar{n}})
\frac{\partial} {\partial \bar{t}_{l}} - e^{-2\bar{t}_m/\bar{n}}
\right] \\
& \left[ \sum^{m}_{l=0} e^{-\bar{t}_l/\bar{n}} \bar{D}^{l,1} \right]
\log\PPn \times \left[ \sum^{m}_{l=0} e^{-\bar{t}_l/\bar{n}}
\bar{D}^{l,1} \right] \left[ \sum^{m}_{l=0} e^{(\bar{t}_l -
\bar{t}_m)/\bar{n}} \bar{D}^{l,1} \right] \left[ \sum^{m}_{l=0}
e^{-\bar{t}_l/\bar{n}} \bar{D}^{l,1} \right] \log\PPn \\
& - \left[ \sum^{m}_{l=0} e^{-2\bar{t}_l/\bar{n}} \bar{D}^{l,2}
+2\bar{n}^{2} \sum^{m}_{l=0} e^{-2\bar{t}_l/\bar{n}} \bar{D}^{l,1} +
\bar{n} \sum^{m}_{l=1} (1 - e^{-2\bar{t}_l/\bar{n}}) \frac{\partial}
{\partial \bar{t}_{l}} - e^{-2\bar{t}_m/\bar{n}} \right] \\
& \left[ \sum^{m}_{l=0} e^{(\bar{t}_l - \bar{t}_m)/\bar{n}}
\bar{D}^{l,1} \right] \log\PPn \times \left[ \sum^{m}_{l=0}
e^{(\bar{t}_l - \bar{t}_m)/\bar{n}} \bar{D}^{l,1} \right] \left[
\sum^{m}_{l=0} e^{-\bar{t}_l/\bar{n}} \bar{D}^{l,1} \right] \left[
\sum^{m}_{l=0} e^{(\bar{t}_l - \bar{t}_m)/\bar{n}} \bar{D}^{l,1}
\right] \log\PPn.
\end{split}
\end{equation}
Since we have commutator formulas
\begin{align}
\left[ \bar{n} (e^{2(\bar{t}_l - \bar{t}_m)/\bar{n}} -
e^{-2\bar{t}_m/\bar{n}}) \frac{\partial}{\partial \bar{t}_l}, \quad
\sum^m_{l=0} e^{-\bar{t}_l/\bar{n}} \bar{D}^{l,1} \right] = &
\sum^m_{l=0} (e^{(\bar{t}_l - 2\bar{t}_m)/\bar{n}} - e^{-(\bar{t}_l
+ 2\bar{t}_m)/\bar{n}})\bar{D}^{l,1}, \\
\left[ \bar{n} \sum^m_{l=1} (1 - e^{-2\bar{t}_l/\bar{n}})
\frac{\partial}{\partial \bar{t}_l}, \quad \sum^m_{l=0}
e^{(\bar{t}_l - \bar{t}_m)/\bar{n}} \bar{D}^{l,1} \right] = &
\sum^m_{l=0} (e^{(\bar{t}_l - 3\bar{t}_m)/\bar{n}} - e^{-(\bar{t}_l
+ \bar{t}_m)/\bar{n}})\bar{D}^{l,1}, \\
\left[ \sum^m_{l=0} e^{-\bar{t}_l/\bar{n}} \bar{D}^{l,1}, \quad
\sum^m_{l=0} e^{2(\bar{t}_l - \bar{t}_m)/\bar{n}} \bar{D}^{l,2}
\right] = & \sum^m_{l=0} e^{(\bar{t}_l - 2\bar{t}_m)/\bar{n}} \bar{D}^{l,1}, \\
\left[ \sum^m_{l=0} e^{(\bar{t}_l - \bar{t}_m)/\bar{n}}
\bar{D}^{l,1}, \quad \sum^m_{l=0} e^{2\bar{t}_l/\bar{n}}
\bar{D}^{l,2} \right] = & \sum^m_{l=0} e^{-(\bar{t}_l +
\bar{t}_m)/\bar{n}} \bar{D}^{l,1},
\end{align}
we can write
(\ref{eq:scaled_equation_for_Dyson_process_before_asymptotics}) as
\begin{equation}
\begin{split}
& \left( \left[ \sum^{m}_{l=0} e^{-\bar{t}_l/\bar{n}} \bar{D}^{l,1}
\right]^2 \left[\sum^{m}_{l=0} e^{2(\bar{t}_l - \bar{t}_m)/\bar{n}}
\bar{D}^{l,2} + 2\bar{n}^{2} \sum^{m}_{l=0} e^{2(\bar{t}_l -
\bar{t}_m)/\bar{n}} \bar{D}^{l,1} + \bar{n} \sum^{m}_{l=1}
(e^{2(\bar{t}_l - \bar{t}_m)/\bar{n}} - e^{-2\bar{t}_m/\bar{n}})
\frac{\partial} {\partial \bar{t}_{l}} \right] \log\PPn \right. \\
& \left. - \left[ \sum^{m}_{l=0} e^{(\bar{t}_l - \bar{t}_m)/\bar{n}}
\bar{D}^{l,1} \right]^2 \left[\sum^{m}_{l=0} e^{-2\bar{t}_l/\bar{n}}
\bar{D}^{l,2} + 2\bar{n}^{2} \sum^{m}_{l=0} e^{-2\bar{t}_l/\bar{n}}
\bar{D}^{l,1} + \bar{n} \sum^{m}_{l=1} (1 - e^{-2\bar{t}_l/\bar{n}})
\frac{\partial} {\partial \bar{t}_{l}} \right] \log\PPn
\right) \\
& \times \left( \left[ \sum^{m}_{l=0} e^{-\bar{t}_l/\bar{n}}
\bar{D}^{l,1} \right] \left[ \sum^{m}_{l=0} e^{(\bar{t}_l -
\bar{t}_m)/\bar{n}} \bar{D}^{l,1} \right] \log\PPn + \bar{n}^2
e^{-\bar{t}_m/\bar{n}} \right) \\
= & \left( \left[ \sum^{m}_{l=0} e^{-\bar{t}_l/\bar{n}}
\bar{D}^{l,1} \right] \left[\sum^{m}_{l=0} e^{2(\bar{t}_l -
\bar{t}_m)/\bar{n}} \bar{D}^{l,2} + 2\bar{n}^{2} \sum^{m}_{l=0}
e^{2(\bar{t}_l - \bar{t}_m)/\bar{n}} \bar{D}^{l,1} + \bar{n}
\sum^{m}_{l=1} (e^{2(\bar{t}_l - \bar{t}_m)/\bar{n}} -
e^{-2\bar{t}_m/\bar{n}}) \frac{\partial} {\partial \bar{t}_{l}}
\right] \log\PPn \right. \\
& \left. - \left[ \sum^{m}_{l=0} e^{(\bar{t}_l -
2\bar{t}_m)/\bar{n}} \bar{D}^{l,1} \right] \log\PPn \right) \times
\left[ \sum^{m}_{l=0} e^{-\bar{t}_l/\bar{n}} \bar{D}^{l,1} \right]^2
\left[ \sum^{m}_{l=0} e^{(\bar{t}_l - \bar{t}_m)/\bar{n}}
\bar{D}^{l,1} \right] \log\PPn \\
& \left( - \left[ \sum^{m}_{l=0} e^{(\bar{t}_l - \bar{t}_m)/\bar{n}}
\bar{D}^{l,1} \right] \left[\sum^{m}_{l=0} e^{-2\bar{t}_l/\bar{n}}
\bar{D}^{l,2} + 2\bar{n}^{2} \sum^{m}_{l=0} e^{-2\bar{t}_l/\bar{n}}
\bar{D}^{l,1} + \bar{n} \sum^{m}_{l=1} (1 - e^{-2\bar{t}_l/\bar{n}})
\frac{\partial} {\partial \bar{t}_{l}} \right] \log\PPn \right. \\
& \left. - \left[ \sum^{m}_{l=0} e^{-(\bar{t}_l +
\bar{t}_m)/\bar{n}} \bar{D}^{l,1} \right] \log\PPn \right) \times
\left[ \sum^{m}_{l=0} e^{-\bar{t}_l/\bar{n}} \bar{D}^{l,1} \right]
\left[ \sum^{m}_{l=0} e^{(\bar{t}_l - \bar{t}_m)/\bar{n}}
\bar{D}^{l,1} \right]^2 \log\PPn.
\end{split}
\end{equation}

Since all terms of the PDE involves $\bar{n}$, we can expand the PDE
with respect to $\bar{n}$, with formulas ($*$ can be $1$ or $2$)
\begin{align}
\sum^{m}_{l=0} e^{-\bar{t}_l/\bar{n}} \bar{D}^{l,1} = &
\sum^{m}_{l=0}\bar{D}^{l,1} - \frac{1}{\bar{n}} \sum^{m}_{l=0}
\bar{t}_l \bar{D}^{l,1} + \frac{1}{2\bar{n}^{2}} \sum^{m}_{l=0}
\bar{t}^2_l \bar{D}^{l,1} - \frac{1}{6\bar{n}^{3}} \sum^{m}_{l=0}
\bar{t}^3_l \bar{D}^{l,1} + O \left( \frac{1}{\bar{n}^{4}} \right), \\
\sum^{m}_{l=0} e^{-2\bar{t}_l/\bar{n}} \bar{D}^{l,*} = &
\sum^{m}_{l=0}\bar{D}^{l,*} - \frac{2}{\bar{n}} \sum^{m}_{l=0}
\bar{t}_l \bar{D}^{l,*} + \frac{2}{\bar{n}^{2}} \sum^{m}_{l=0}
\bar{t}^2_l \bar{D}^{l,*} - \frac{4}{3\bar{n}^{3}} \sum^{m}_{l=0}
\bar{t}^3_l \bar{D}^{l,*} + O \left( \frac{1}{\bar{n}^{4}} \right),  \\
\sum^{m}_{l=0} e^{(\bar{t}_l - \bar{t}_m)/\bar{n}} \bar{D}^{l,1} = &
\sum^{m}_{l=0}\bar{D}^{l,1} + \frac{1}{\bar{n}} \sum^{m}_{l=0}
(\bar{t}_l - \bar{t}_m) \bar{D}^{l,1} + \frac{1}{2\bar{n}^{2}}
\sum^{m}_{l=0} (\bar{t}_l - \bar{t}_m)^{2} \bar{D}^{l,1} \notag \\
& + \frac{1}{6\bar{n}^{3}} \sum^{m}_{l=0} (\bar{t}_l -
\bar{t}_m)^{3} \bar{D}^{l,1} + O \left( \frac{1}{\bar{n}^{4}} \right), \\
\sum^{m}_{l=0} e^{2(\bar{t}_l - \bar{t}_m)/\bar{n}} \bar{D}^{l,*} =
& \sum^{m}_{l=0}\bar{D}^{l,*} + \frac{2}{\bar{n}} \sum^{m}_{l=0}
(\bar{t}_l - \bar{t}_m) \bar{D}^{l,*} + \frac{2}{\bar{n}^{2}}
\sum^{m}_{l=0} (\bar{t}_l - \bar{t}_m)^{2}
\bar{D}^{l,*} \notag \\
& + \frac{4}{3\bar{n}^{3}} \sum^{m}_{l=0} (\bar{t}_l -
\bar{t}_m)^{3} \bar{D}^{l,*} + O \left( \frac{1}{\bar{n}^{4}}
\right),
\end{align}
\begin{align}
\bar{n} \sum^m_{l=1} (e^{2(\bar{t}_l - \bar{t}_m)/\bar{n}} -
e^{2\bar{t}_m/\bar{n}}) \frac{\partial}{\partial \bar{t}_l} = & 2
\sum^m_{l=1} \bar{t}_l \frac{\partial}{\partial \bar{t}_l} +
\frac{2}{\bar{n}} \sum^m_{l=1} \bar{t}_l (\bar{t}_l - 2\bar{t}_m)
\frac{\partial}{\partial \bar{t}_l} + O \left( \frac{1}{\bar{n}^{2}} \right), \\
\bar{n} \sum^m_{l=1} (1 - e^{-2\bar{t}_l/\bar{n}})
\frac{\partial}{\partial \bar{t}_l} = & 2 \sum^m_{l=1} \bar{t}_l
\frac{\partial}{\partial \bar{t}_l} - \frac{2}{\bar{n}} \sum^m_{l=1}
\bar{t}^2_l \frac{\partial}{\partial \bar{t}_l} + O \left(
\frac{1}{\bar{n}^{2}} \right).
\end{align}
Although the left hand side of
(\ref{eq:scaled_equation_for_Dyson_process_before_asymptotics})
contains $O(\bar{n}^4)$ terms, after careful calculation we find all
$O(\bar{n}^4)$, $O(\bar{n}^3)$ and $O(\bar{n}^2)$ terms disappear,
and the equation becomes
\begin{equation} \label{eq:limit_of_equation_of_Airy_process}
\begin{split}
& \left[ \sum^m_{l=0} \bar{D}^{l,1} \right]^2 \left[ \sum^m_{l=0}
(\bar{t}_m - 2\bar{t}_l) \bar{D}^{l,2} + \sum^m_{l=0} ((\bar{t}_m -
\bar{t}_l)^3 - \bar{t}^3_l)^3 \bar{D}^{l,1} + 2\sum^m_{l=1}
\bar{t}_l (\bar{t}_m - \bar{t}_l) \frac{\partial}{\partial
\bar{t}_l} \right] \log \PPn \notag \\
& + \left[ \sum^m_{l=0} \bar{D}^{l,1} \right] \left[ \sum^m_{l=0}
(2\bar{t}_l - \bar{t}_m) \bar{D}^{l,1} \right] \left[ \sum^m_{l=0}
\bar{D}^{l,2} + \sum^m_{l=0} (\bar{t}^2_l + (\bar{t}_m -
\bar{t}_l)^2) \bar{D}^{l,1} + 2\sum^m_{l=1} \bar{t}_l
\frac{\partial}{\partial \bar{t}_l} \right] \log \PPn \notag \\
& + \left[ \sum^m_{l=0} \bar{D}^{l,1} \right] \left[ \left[
\sum^m_{l=0} \bar{t}_l \bar{D}^{l,1} \right] \left[ \sum^m_{l=0}
(\bar{t}_m - \bar{t}_l)^2 \bar{D}^{l,1} \right] - \left[
\sum^m_{l=0} (\bar{t}_m - \bar{t}_l) \bar{D}^{l,1} \right] \left[
\sum^m_{l=0} \bar{t}^2_l \bar{D}^{l,1} \right]
\vphantom{\left[ \sum^m_{l=0} \bar{D}^{l,1} \right]^2}
\right] \log \PPn \notag \\
& + 2 \left[ \sum^m_{l=0} \bar{t}_l \bar{D}^{l,1} \right] \left[
\sum^m_{l=0} (\bar{t}_m - \bar{t}_l) \bar{D}^{l,1} \right] \left[
\sum^m_{l=0} (2\bar{t}_l - \bar{t}_m) \bar{D}^{l,1} \right]
\log \PPn \notag \\
= & \left\{ \left[ \sum^m_{l=0} (2\bar{t}_l - \bar{t}_m)
\bar{D}^{l,1} \right] \left[ \sum^m_{l=0} \bar{D}^{l,1} \right] \log
\PPn, \left[ \sum^m_{l=0} \bar{D}^{l,1} \right] \left[ \sum^m_{l=0}
\bar{D}^{l,1} \right] \log \PPn \right\}_{\sum^m_{l=0}
\bar{D}^{l,1}} + O \left( \frac{1}{\bar{n}} \right).
\end{split}
\end{equation}

The term $O \left( \frac{1}{\bar{n}} \right)$ in
(\ref{eq:limit_of_equation_of_Airy_process}) is a quadratic function
in term of $\log \PPn$ and its derivatives with coefficients $O
\left( \frac{1}{\bar{n}} \right)$. By the definition of $\PPAiry$ in
(\ref{eq:definition of Airy process}) and the convergence result in
\cite{Adler-van_Moerbeke05}, we take the limit $n \rightarrow
\infty$, and get the PDE (\ref{eq:equation_for_Airy_process}) after
the changing of notations, i.e., cleaning all ``bars'' for variables
and operators.

\bibliographystyle{plain}
\bibliography{bibliography}
\end{document}